\documentclass[12pt, reqno]{amsart}
\usepackage{amsmath, amsthm, amscd, amsfonts, amssymb, graphicx, color}
\usepackage[bookmarksnumbered, plainpages]{hyperref}
\usepackage{mathtools}
\usepackage{scalerel}
\usepackage{nccmath}
\usepackage{bm}
\usepackage{commath}
\usepackage{relsize}

\newtheoremstyle{break}
{\topsep}{\topsep}%
{\itshape}{}%
{\bfseries}{}%
{\newline}{}%

\theoremstyle{break}
\newtheorem{definition}{Definition}[section]
\newtheorem{theorem}[definition]{Theorem}
\newtheorem{lemma}[definition]{Lemma}
\newtheorem{prop}[definition]{Proposition}
\newtheorem{corollary}[definition]{Corollary}

\newtheorem{thmx}{Theorem}
\newtheorem{defx}[thmx]{Definition}
\newtheorem{lemx}[thmx]{Lemma}

\newcommand\blfootnote[1]{%
	\begingroup
	\renewcommand\thefootnote{}\footnote{#1}%
	\addtocounter{footnote}{-1}%
	\endgroup
}

\newenvironment{ProofOf}{$ $}{\hfill$\square$}

\newcommand{\D}{\mathbb{D}}
\newcommand{\C}{\mathbb{C}}
\newcommand{\F}{\mathcal{F}}
\newcommand{\N}{\mathbb{N}}
\newcommand{\R}{\mathbb{R}}

\newcommand{\e}{\text{e}}
\newcommand{\fn}{(f_n)_n}
\newcommand{\HD}{\mathcal{H}(\D)}
\newcommand{\MD}{\mathcal{M}(\D)}
\newcommand{\HkM}{\mathcal{H}_{k,M}}
\newcommand{\HkMD}{\mathcal{H}_{k,M}'}
\newcommand{\HkML}{\mathcal{H}_{k,M}^{f'/f}}
\newcommand{\HkMlog}{\mathcal{H}_{k,M}^{f''/f'}}
\newcommand{\MkM}{\mathcal{M}_{k,M}}
\newcommand{\MkMD}{\mathcal{M}_{k,M}'}
\newcommand{\MkML}{\mathcal{M}_{k,M}^{f'/f}}
\newcommand{\MkMlog}{\mathcal{M}_{k,M}^{f''/f'}}

\begin{document}
	\setcounter{page}{1}
	\mathtoolsset{showonlyrefs,showmanualtags}
	
	\title{Generalized Schwarzians and Normal Families}	
	\author{Matthias Grätsch}
	\maketitle	
	\begin{abstract}
		We study families of analytic and meromorphic functions with bounded generalized Schwarzian derivative $S_k(f)$. We show that these families are quasi-normal.  Further, we investigate associated families, such as those formed by derivatives and logarithmic derivatives, and prove several \mbox{(quasi-)}normality results. Moreover, we derive a new formula for $S_k(f)$, which yields a result for families $\F\subseteq\HD$ of locally univalent functions that satisfy
		\begin{align}
			S_k(f)(z)\neq b(z)\qquad \text{for some }b\in\mathcal{M}(\D)\text{ and all } f\in\F,\,z\in\C
		\end{align}
		and for entire functions $f$ with $S_k(f)(z)\neq0$ and $S_k(f)(z)\neq\infty$ for all $z\in\C$.\\
		The classical Schwarzian derivative $S_f$ is contained as the case $k=2$.
	\end{abstract} 
	\section{Introduction and Main Results}\label{se:Main}
	\noindent 
	Throughout this paper, we denote the set of all holomorphic functions on a domain $D\subseteq\C$ by $\mathcal{H}(D)$. Likewise, we write $\mathcal{M}(D)$ for the set of all meromorphic functions on $D$. Further, we denote the open unit disk by ${\D\coloneq\{z\in\C\,\,:\,\,\abs{z}<1\}}$. Moreover, when referring to zeros or poles, we use ``(CM)'' to indicate that multiplicity is counted, and ``(IM)'' when it is ignored.\\
	\noindent
	Let $D\subseteq\C$ be a domain. A family $\mathcal{F}\subseteq\mathcal{M}(D)$ is said to be \textbf{quasi-normal on} $D$, if for every sequence $\fn\subseteq\mathcal{F}$, there exists a subsequence $(f_{n_k})_{k}\subseteq\fn$ and an exceptional set $E\subseteq D$ with no accumulation point in $D$, such that $(f_{n_k})_{k}$ converges locally uniformly in $D\backslash E$ (with respect to the spherical metric).\\
	A family $\mathcal{F}\subseteq\MD$ is said to be \textbf{quasi-normal at} $z_0\in D$ if there exists a neighborhood $U\subseteq D$ around $z_0$, such that the restricted family $\{f|_U\,:\, f\in\mathcal{F}\}$ is quasi-normal.\medskip
	
	\noindent
	The theory of normal families offers a multitude of different criteria to check if a given family $\F\subseteq\MD$ is normal. Probably the most notable one is Marty's theorem, which states that a family $\F$ is normal, if and only if the family of spherical derivatives $\big\{\frac{\vert f'\vert}{1+\vert f \vert^2}\,:\,f\in\F\big\}$ is locally uniformly bounded. This draws a connection between the normality of a family and a particular differential inequality. Similarly, there are several (more or less) related results, which connect (quasi-)normality to other differential operators and inequalities, such as\blfootnote{$\!\!\!\!\!\!\!$2020 Mathematics Subject Classification: Primary 30D45; Secondary 30D30, 30C45\\
		Key words: normal family, Schwarzian derivative, quasi-normality, value distribution theory}
	\begin{align}
		\frac{|f^{(k)}|}{1+|f^{(j)}|^\alpha}\leq M,\qquad  M\leq\frac{|f^{(k)}|}{1+|f^{(j)}|^\alpha}\qquad\text{ or }\qquad f^n(z)+f^{(k)}(z)\neq 0
	\end{align}
	for suitable choices of $M,\alpha, k,j,n$ and all admissible $z$ (see for example \cite{GrahlNevo}, \cite{GrahlManketNevo}, \cite{BarGrahlNevo},  \cite{NevoShemTov} and \cite{ChenHua}).\pagebreak[3]
	
	\noindent
	This paper seizes this idea and investigates (quasi-)normality in relation to inequalities involving a generalization of the Schwarzian derivative by M. Chuaqui, J. Gröhn and J. Rättyä (see \cite{ChuaquiGrönRättyä}). This generalization is anchored around the relationship between the Schwarzian derivative and the so called Schwarzian differential equation $y''+p_0\,y=0$. Other generalizations of the Schwarzian derivative, which focus on different aspects of the Schwarzian derivative, can be found in \cite{KimSugawa},\cite{Schippers} or \cite{Tamanoi}.
	\medskip
	\begin{defx}[\mbox{\cite[p. 340]{ChuaquiGrönRättyä}}]\label{de:GeneralizedSchwarzian}
		Let $f\in\mathcal{M}(\D)$ be non-constant. For $n\in\N$ and $k\in\N\backslash\{1\}$, we define:
		\begin{align*}
			S_{2,n}(f)\coloneq\frac{f''}{f'} \qquad \text{ and } \qquad 	S_{k+1,n}(f)\coloneq\Big(S_{k,n}(f)\Big)'-\frac{1}{n}\frac{f''}{f'}S_{k,n}(f).
		\end{align*} 
		Now, for all $k\in\N$, we call 
		\begin{align*}
			S_k(f)\coloneq S_{k+1,k}(f)
		\end{align*}
		the generalized Schwarzian derivative of order k.
	\end{defx}
	\medskip
	\noindent
	Note that the classical Schwarzian derivative $S_f$ is contained as the case $k=2$, since
	\begin{align*}
		S_2(f) = S_{3,2}(f) = \left(S_{2,2}(f)\right)'-\frac{1}{2}\frac{f''}{f'}S_{2,2}(f)= \bigg(\frac{f''}{f'}\bigg)'-\frac{1}{2}\bigg(\frac{f''}{f'}\bigg)^2 = S_f.
	\end{align*}
	For constant functions $c\in\HD$, we define $S_k(c)\equiv\infty$ for all $k\in\N$. This differs from \cite{ChuaquiGrönRättyä}, where constant functions were not mentioned at all. However, our version allows us to circumvent degenerate cases in the future. It is further motivated by the observation that $S_{k}(f)$ exhibits a $(k-1)$-fold pole whenever the derivative of a non-constant $f\in\MD$ vanishes. We will give a short proof of this result in Section \ref{se:Lemmas}, Proposition \ref{pr:Poles}.\\
	The connection between this generalization of the Schwarzian derivative and the Schwarzian differential equation becomes apparent in the following theorem.\medskip
	\begin{thmx}[\mbox{\cite[Lemma 3 and Lemma 5]{ChuaquiGrönRättyä}}]\label{th:ChuaquiGrönRättyä}
		Let $f\in\MD$ and $k\in\N$. Then the following conditions are equivalent:
		\begin{enumerate}
			\item[(a)] $S_k(f)\in\HD$.
			\item[(b)] $f'=1/h^k$, for some $h\in\HD$, which satisfies the differential equation ${y^{(k)}+p_0\, y=0}$ for some $p_0\in\HD$.
			\item[(c)] $f'=1/h^k$ for some $h\in\HD$, where every zero $z_0\in\D$ of $h$ with multiplicity $m$ is also a zero of $h^{(k)}$ of multiplicity at least $m$.
		\end{enumerate}
		If condition \textnormal{(b)} is true, then we can specify $p_0=S_k(f)/k$ and ${S_k(f)=-k\,h^{(k)}/h}$.
	\end{thmx}\medskip\noindent
	Again, this result differs slightly from the one given in \cite[Lemma 5]{ChuaquiGrönRättyä}. There, it is additionally required that $f'$ is non-vanishing. However, with our exclusion of constant functions and Proposition \ref{pr:Poles}, we can drop this assumption.\bigskip
	
	\noindent
	To abbreviate our notation, we define for each $M\in\R$ and $k\in\N$:
	\begin{align}
		\MkM\coloneq\{f\in\MD\;:\,\Vert S_k(f)\Vert_\infty\leq M\}\qquad\text{and}\qquad \HkM\coloneq \MkM\cap\HD,
	\end{align} 
	where $\Vert\cdot\Vert_\infty$ denotes the supremum norm on $\D$. Since no constant function is contained in $\HkM$ or $\MkM$, we can also consider the following families:
	\begin{align}
		&\mathcal{M}_{k,M}'\coloneq\{f'\,:\, f\in\mathcal{M}_{k,M}\}&&\text{ and }&&\,\mathcal{H}_{k,M}'\coloneq\{f'\,:\,  	f\in\mathcal{H}_{k,M}\},\\
		&\mathcal{M}_{k,M}^{f'/f}\coloneq\{g'/g\,:\,  g\in\mathcal{M}_{k,M}\} &&\text{ and 	}&&\,\mathcal{H}_{k,M}^{f'/f}\coloneq\{g'/g\,:\,  g\in\mathcal{H}_{k,M}\},\\
		&\mathcal{M}_{k,M}^{f''/f'}\coloneq\{g''/g'\,:\,  g\in\mathcal{M}_{k,M}\} &&\text{ and 	}&&\,\mathcal{H}_{k,M}^{f''/f'}\coloneq\{g''/g'\,:\,  g\in\mathcal{H}_{k,M}\}.
	\end{align}
	In the case of the classical Schwarzian, the quasi-normality of $\mathcal{H}_{2,M}$ can be shown for any $M\in\R$ as a consequence of \cite[Theorem 1.1(b)]{MaMeijaMinda}. The quasi-normality of $\mathcal{M}_{2,M}$ was later proven in \cite[Theorem 1.4]{Grätsch} for all $M\in\R$. There, it is also shown for any $M\in\R$ that $\mathcal{M}_{2,M}^{f''/f'}$ is normal and $\mathcal{M}_{2,M}^{f'/f'}$ is quasi-normal.\\
	In this paper, we will generalize these results for all $k\in\N$. To achieve this, we consider families of holomorphic functions first and obtain the following theorem.\medskip
	\begin{theorem}\label{th:HkM}
		The following statements hold for all $k\in\N$ and $M\in\R$:
		\begin{itemize}
			\item[] $\HkMlog$ is locally uniformly bounded.
			\item[] $\HkML$ is normal and no sequence in $\HkML$ converges to $\infty$.
			\item[] $\HkMD$ is normal.
			\item[] $\HkM$ is quasi-normal.
		\end{itemize}
	\end{theorem}
	\medskip\noindent
	Later, in Lemma \ref{le:FinitePoles}, we will see that the number of poles of each $f\in\MkM$ is bounded by a constant that depends only on $k$ and $M$. Thus, every sequence $\fn\subseteq\MkM$ has a subsequence such that the poles of all its functions form a set whose accumulation points are isolated. Since (quasi-)normality is a local property, we are able to treat $\MkM$ mostly like its holomorphic subset $\HkM$ and prove the following theorem:\medskip
	
	\begin{theorem}\label{th:MkM}
		The following statements hold for all $k\in\N$ and $M\in\R$:
		\begin{itemize}
			\item[] $\MkMlog$ is quasi-normal and no sequence in $\MkMlog$ converges to $\infty$.
			\item[] $\MkML$ is quasi-normal and no sequence in $\MkML$ converges to $\infty$.
			\item[] $\MkMD$ is quasi-normal.
			\item[] $\MkM$ is quasi-normal.
		\end{itemize}
	\end{theorem}
	\medskip
	\noindent
	\cite[Lemma 4]{ChuaquiGrönRättyä} shows that $S_k(f)\equiv 0$ if and only if $f'=1/p^k$, where $p$ is a polynomial with $\deg p\leq k-1$. Thus, $S_k(f_n)\equiv 0$ for $f_n(z)=nz$, for all $k\in\N$, so neither $\HkM$ nor $\MkM$ are normal for $k\in\N$ and $M\in\R_0^+$. \\
	Similarly, $\MkMD$ is not normal for $k\in\N\backslash\{1\}$ and $M\in\R_0^+$, since $S_k(g_n)\equiv0$ for 
	\begin{align}
		g_n(z)\coloneq \frac{-1}{n^k(k-1)z^{k-1}}\qquad\text{with} \qquad g_n'(z)=\frac{1}{(nz)^k}.
	\end{align}
	For $k\in\N\backslash\{1\}$ and $M\in\R_0^+$, the family $\MkML$ is not normal either. To see this, we consider the sequence
	\begin{align}
		h_n(z)\coloneq \frac{1}{z^{k-1}}+n\quad\text{with}\quad h_n'(z)=\frac{1-k}{z^k}\quad\text{and}\quad \frac{h_n'}{h_n}(z)=\frac{1-k}{z(1+nz^{k-1})}.
	\end{align}
	It is unknown to the author, whether $\MkMlog$ is normal. However, by using \mbox{Theorem \ref{th:HkM}}, it is possible to show that $\MkMlog$ is normal at $z_0\in\D$ if and only if there exists a neighborhood of $z_0$, where each $f\in\MkM$ has at most one \mbox{pole (IM)}.  Theorem \ref{th:Disconjucacy} and the proof of \cite[Theorem 1.4]{Grätsch} show that these conditions hold for $k=2$, but it remains open whether this is true for $k\geq3$.\\~\\
	For our final results, we will regard $S_k(f)$ as a differential polynomial by showing the following formula for the generalized Schwarzian derivative.\medskip
	\begin{lemma}\label{le:SchwarzianDifferentialPolynomial}
		Let $f\in\MD$ be a non-constant, meromorphic function and $k\in\N$. Then
		\begin{align}\label{eq:FormulaS_k(f)}
			S_k(f)=\mathlarger{\sum}_{\mathclap{(n_1,\hdots,n_k)\,\in\,\Lambda}}\,\,\,\,\,\,\Bigg(\frac{-k\cdot 	k!}{\prod_{j=1}^k (-k\cdot j!)^{n_j}\cdot n_j!}\cdot\prod\limits_{j=1}^k\Bigg(\bigg(\frac{f''}{f'}\bigg)^{(j-1)}\Bigg)^{n_j}\Bigg),
		\end{align}
		where $\Lambda$ is the set of all tuples $(n_1,\hdots,n_k)$ (with $n_r\in\N_0$ for all $r=1,\hdots,k$) that satisfy $\sum_{r=1}^{k}r\cdot n_r=k$.
	\end{lemma}\medskip\noindent\sloppy
	For $k\in\N\backslash\{1\}$, we can extract the summands where ${(n_1\hdots, n_k)=(k,0,\hdots, 0)}$ and ${(n_1,\hdots,n_k)=(0,\hdots,0,1)}$, to obtain that 
	\begin{align}\label{eq:DifferentialPolynomialP}
		S_k(f)
		=\frac{(-1)^{k+1}}{k^{k-1}}\cdot g^k+g^{(k-1)}+P[g]
	\end{align}
	for $g\coloneq f''/f'$ and some differential polynomial $P$. \\
	This form is reminiscent of the well known condition ${af^n(z)+f^{(m)}(z)\neq b}$ for $a,b\in\C$, $a\neq0$ and large enough $n,m\in\N$ with $n>m$. \mbox{W. K. Hayman} was the first to study this condition for $f\in\mathcal{M}(\C)$ (see \cite{Hayman}). In the spirit of Bloch's principle, his results on value distribution were later extended to normal families: in the analytic case by D. Drasin in \cite{Drasin}, and in the meromorphic case by J. K. Langley in \cite{Langley}. Subsequently, H. Chen and X. Hua demonstrated in \cite{ChenHua} that the constant value $b$ can be replaced by an exceptional function $B\in\mathcal{M}(\D)$. This result was further extended to more general differential polynomials (see Theorem~\ref{th:GrahlNormal} below) by J. Grahl in \cite{Grahl}. Now, we will show that the differential polynomial $P$ from equation~\eqref{eq:DifferentialPolynomialP} fulfills the exact conditions of Theorem~\ref{th:GrahlNormal}. This leads to the following result:\medskip 
	\begin{theorem}\label{th:Exeption}
		Let $\F\subseteq\HD$ be a family of locally univalent functions, $k\in\N\backslash\{1\}$ and ${b\in\MD}$ be a meromorphic function, such that every $f\in\F$ satisfies
		\begin{align}\label{eq:Exception}
			S_k(f)(z)\neq b(z) \qquad\text{for all}\qquad z\in\D.
		\end{align}
		Then $\F''/\F'\coloneq\{f''/f'\;:\,f\in\F\}$ is a normal family.
	\end{theorem}\medskip\noindent
	Note that $\F''/\F'$, unlike $\HkMlog$, does not have to be locally uniformly bounded. Therefore, the family $\F$ does not need to be quasi-normal, as demonstrated by the family $\{z\mapsto e^{nz}\,:\,n\in\N\}$.\pagebreak
	
	\noindent
	However, if - in addition to \eqref{eq:Exception} - we know that $\F''/\F'$ is pointwise bounded in a single point, then it follows that $\F''/\F'$ is locally uniformly bounded. Now, using formula \eqref{eq:FormulaS_k(f)}, we can see that the family $S_k(\F)\coloneq\{S_k(f)\,:\,f\in\F\}$ is locally uniformly bounded as well. This, in turn, allows us to apply Theorem \ref{th:HkM} to $\F$. Thus, we can regard this as a ``self-improving result'', where we require that:
	\begin{itemize}
		\item $S_k(\F)$ omits a function.
		\item $\F''/\F'$ is bounded in a single point.
	\end{itemize}
	and obtain that:
	\begin{itemize}
		\item $S_k(\F)$ is locally uniformly bounded.
		\item $\F''/\F'$ is locally uniformly bounded.
		\item $\F'/\F$ and $\F'$ are normal, while $\F$ is quasi-normal.
	\end{itemize}
	\medskip
	We should note that Theorem \ref{th:Exeption} can not be extended to families of meromorphic functions. To see this, consider  $f_n(z)\coloneq((2z)^n-1)^{-1}$ for $n\in\N$ and $z\in\D$ with
	\begin{align}
		f_n'(z)=-\frac{n\,2^n\,z^{n-1}}{\big((2z)^n-1\big)^2}\qquad\text{and}\qquad 
		f_n''(z)=\frac{n\,2^n\,z^{n-2}\big(n-1+(n+1)(2z)^n\big)}{\big((2z)^n-1\big)^3}.
	\end{align}
	Clearly, each $f_n$ has a pole of order $1$ in $z_0=1/2$ and is locally injective in $\D\backslash\{0\}$. Additionally, for $z_n\coloneq \frac{1}{2}\sqrt[n]{\frac{n-1}{n+1}}\,e^{i\pi/n}$, we have $f_n''(z_n)=0$, while $(z_n)_n$ converges to $1/2$. Hence, the sequence $(f_n''/f_n')_n$ is not normal at $1/2$.\\
	On the other hand, we can use the fact that the classical Schwarzian derivative is invariant under Möbius transformations, and calculate 
	\begin{align}
		S_2(f_n)(z)
		=S_2(z^n)
		=\frac{1-n^2}{2\,z^2}.
	\end{align}
	Thus, $\big(S_k(f_n)\big)_{n\geq2}$ omits the value $0$, which shows that we can not extend \mbox{Theorem \ref{th:Exeption}} to families of meromorphic functions.\medskip
	
	\noindent
	Still, there is a corresponding value distribution result for entire functions. Similarly to Theorem \ref{th:Exeption}, this corollary relies heavily on the results from \cite{Grahl}.\medskip
	\begin{corollary}\label{co:ExceptionalEnitre}
		Let $f\in\mathcal{H}(\C)$ be an entire function with
		\begin{align}
			S_k(f)(z)\neq0\quad\text{ and }\quad S_k(f)(z)\neq\infty\qquad\text{for all } z\in\C.
		\end{align}
		Then there are $a,b,c\in\C$ with $a,b\neq0$ and $f(z)=ae^{bz}+c$.
	\end{corollary}\medskip\noindent
	Based on the counterexample given in \cite[p. 34]{Hayman}, we consider the locally univalent function ${f(z)=\exp(\exp(cz)/c)}$ with $c\in\C\backslash\{0\}$. A straightforward computation shows that its classical Schwarzian derivative is
	\begin{align}
		S_2(f)(z)=-\frac{e^{2cz}}{2}-\frac{c^2}{2},
	\end{align}
	so the exceptional value $0$ cannot be replaced by another value.\pagebreak

	\noindent
	Likewise, Corollary \ref{co:ExceptionalEnitre} does not extend to meromorphic functions. For $k=2$, this is already evident from compositions of Möbius transformations with exponential functions.	Now, a natural question is whether any ${f\in\mathcal{M}(\C)}$ with $S_k(f)(z)\neq0$ and $S_k(f)(z)\neq\infty$ for all $z\in\C$ must have the form ${f(z)=T(ae^{bz}+c)}$ for some invariant function $T$ of the $k$-th generalized Schwarzian derivative. However, Proposition \ref{pr:Counterexample} provides a counterexample by showing that such a generalization does not hold for $k=2$.

	\section{Auxiliary Lemmas and Results}\label{se:Lemmas}
	\noindent
	First, we show that non-constant functions with an analytic generalized Schwarzian derivative have non-vanishing derivatives.\medskip
	\begin{prop}\label{pr:Poles}
		If $f\in\MD$ is non-constant and $f'$ has a zero in $z_0\in\D$, then $S_k(f)$ has a $k$-fold pole in $z_0$. 
	\end{prop}\medskip
	\begin{proof}
		Suppose that $f'$ has an $m$-fold zero in $z_0$. Then $S_{2,n}(f)=f''/f'$ has a simple pole in $z_0$ for all $n\in\N$. Next, we assume that $S_{k,n}(f)$ has a $(k-1)$-fold pole in $z_0$ for $k\geq3$, i.e. 
		\begin{align}
			S_{k,n}(f)(z)=\mfrac{h(z)}{(z-z_0)^{k-1}} \quad\text{for some analytic }h \text{ in a neighborhood } U\text{ of } z_0.
		\end{align}
		Then we can calculate\vspace{-5pt}
		\begin{align*}
			S_{k+1,n}(f)(z)
			=\frac{1}{(z-z_0)^k}\bigg(\Big(1-k-\frac{m}{n}\Big)h(z_0)+\hdots\bigg)\qquad\text{for } z\in U,
		\end{align*}
		and conclude inductively that $S_{k,n}(f)$ has a $k$-fold pole for all ${k,n\in\N\backslash\{1\}}$.
	\end{proof}\medskip
	\noindent
	The proof of Theorem \ref{th:HkM} relies heavily on the following theorem by W. Schwick.\medskip
	\begin{thmx}[\mbox{\cite[Theorem 5.4]{Schwick}}]\label{th:Schwick}
		Let $(h_n)_n\subseteq\HD$ be a sequence of non-vanishing functions, and let $k\in\N$. If  $(h^{(k)}_n/h_n)_n$ converges locally uniformly to some $\psi\in\mathcal{H}(\D)$, then $(h_n'/h_n)_n$ is locally uniformly bounded in $\D$.
	\end{thmx}\medskip\noindent
	W. Schwick initially stated this result in a slightly weaker form, establishing only that $(h_n'/h_n)_n$ is normal. However, his proof shows that the Nevanlinna characteristic of $(h_n'/h_n)_n$ is locally uniformly bounded, which, as shown in \mbox{\cite[Theorem 1.13]{Schwick}}, implies that the sequence itself is locally uniformly bounded.\medskip
		
	\noindent Theorem \ref{th:Schwick} will be used to show that $\MkMlog$ is locally uniformly bounded for all $k\in\N$ and $M\in\R$. Then the following result allows us to transfer the convergence properties of sequences in $\HkMlog$ to respective subsequences in $\HkML$ and $\HkM$.\medskip
	
	\begin{lemx}[\mbox{\cite[Lemma 2.4]{Grätsch}}]\label{le:Lift}
		Let $E\subseteq\D$ be a set without an accumulation point in $\D$ and $\fn\subseteq\MD$ with:
		\begin{enumerate}
			\item[(1.)] $(f_n''/f_n')_n$ converges locally uniformly on $\D\backslash E$ to some $\psi\in\mathcal{H}(\D\backslash E)$.
			\item[(2.)] $f_n'$ is zero-free for all $n\in\N$.
		\end{enumerate}
		Then $\fn$ and $(f_n'/f_n)_n$ are quasi-normal on $\D$, and no subsequence of $(f_n'/f_n)_n$ converges to $\infty$.	Moreover, if $E=\emptyset$, $(f_n'/f_n)_n$ is normal.
	\end{lemx}\medskip
	
	\noindent
	Next, we will estimate the maximal number of poles of the functions in $\MkM$ by using the differential equation in Theorem \ref{th:ChuaquiGrönRättyä}(b). Here, the concept of \textit{disconjugate differential equations} will be useful.\medskip
	\begin{defx}
		Let $D\subseteq\C$ be a domain, $k\in\N$ and $p_0,\hdots,p_{k-1}\in\mathcal{H}(D)$. We say that
		\begin{align}\label{eq:ODE}
			y^{(k)}+p_{k-1}\cdot y^{(k-1)}+p_{k-2}\cdot y^{(k-2)}+\hdots+p_{0}\cdot y=0
		\end{align}
		is disconjugate in D, if no non-trivial solution has more than $k-1$ zeros (CM).
	\end{defx}\medskip
	\noindent 
	A classical result concerning the Schwarzian derivative states that if $f_1,f_2\in\HD$ are linearly independent solutions of the differential equation $y''+p_0\cdot y=0$, then ${f\coloneq f_1/f_2}$ satisfies $S_2(f)=p_0/2$ (cf. \cite[Theorem~6.1]{Laine}). As a consequence, any linear combination $g=c_1f_1+c_2f_2$ vanishes at a point $z_0$ if and only if ${f(z_0)=-c_2/c_1}$. Therefore, $f$ attains some value $n$ times if and only if there exists a non-trivial solution of $y''+p_0\cdot y=0$ with $n$ zeros.\\
	This observation is due to Z. Nehari and his paper on the Schwarzian derivative and univalence (see \cite[p. 546]{NehariUnivalence}). Since then, numerous results provided alternative criteria for disconjugacy. One such result is the following theorem. However, with minor modifications, Theorem~\ref{th:Disconjucacy} and  Lemma~\ref{le:FinitePoles} could also be derived from \cite[Theorem~2]{Hadass} or \cite[p. 328]{Nehari}.
	\medskip
	\begin{thmx}[\mbox{\cite[p. 723]{Kim}}]\label{th:Disconjucacy}
		Let $D\subseteq\C$ be a convex domain with $\delta\coloneq\text{diam}(D)<\infty$, $k\in\N$ and let $p_0\in\mathcal{H}(D)$ be a holomorphic function with 
		\begin{align}\label{k!/delta^k}
			|p_0(z)|<\frac{k!}{\delta^k}\qquad\text{ for all } \qquad z\in D.
		\end{align} 
		Then the differential equation $y^{(k)}(z)+p_0(z)y(z)=0$ is disconjugate in $D$.
	\end{thmx}\medskip
	\noindent
	Based on this result, we can show the following lemma.\medskip
	
	\begin{lemma}\label{le:FinitePoles}
		For $k\in\N$ and $M\in\R$, there exists $N\in\N$, such that all $f\in\MkM$ and all $f'\in\MkMD$ have at most $N$ poles (CM).
	\end{lemma}\medskip
	
	\begin{proof}
		For $f\in\MkM$, it suffices to show that $f'$ has a limited number of poles.
		The case $M=0$ was already established in Theorem \cite[Lemma~4]{ChuaquiGrönRättyä}, so we assume $M>0$. Due to Theorem~\ref{th:ChuaquiGrönRättyä}, we can write $f'=1/h^k$, where $h\not\equiv0$ is a solution of the differential equation
		\begin{align}\label{eq:SchwarzianDisconjugate}
			y^{(k)}(z)+\frac{S_k(f)(z)}{k}\cdot y(z)=0.
		\end{align}
		Now, we choose $\delta>0$ with $\delta<\Big(\mfrac{k\cdot k!}{M}\Big)^{1/k}$. This implies ${\mfrac{\Vert S_k(f)\Vert_\infty}{k}\leq\mfrac{M}{k}<\mfrac{k!}{\delta^k}}$, 
		so Theorem \ref{th:Disconjucacy} shows that \eqref{eq:SchwarzianDisconjugate} is disconjugate in every convex set $C_\delta\subseteq\D$ with diameter $\delta$. Consequentially, $h$ can have at most $k-1$ zeros in $C_\delta$. By covering the unit disk $\D$ with $\tilde{N}$ of such convex sets, we obtain the bound $N\coloneq \tilde{N}(k-1)$ for the number of zeros of $h$ (CM).
	\end{proof}
	\pagebreak
	\noindent
	For the reader’s convenience, we also restate the main results of \cite{Grahl} with modified notation, to prevent conflicts with the notation used in this paper.\medskip
	
	\begin{thmx}[\mbox{\cite[Theorem~3]{Grahl}}]\label{th:GrahlNormal}
		Let $\F\subseteq\HD$, $\ell\in\N$ and $k\in\N\backslash\{1\}$. Let $a,b,a_1,\hdots,a_N\in\MD$ with $a\not\equiv0$, and suppose that all poles of $a$ have multiplicity at most $k-1$.\\
		Consider a differential polynomial of the form
		\begin{align}
			P[u]=\sum\limits_{\mu=1}^N a_\mu\cdot\prod\limits_{j=1}^{s_\mu} u^{(\omega_{\mu,j})},
		\end{align}
		where $s_\mu, \omega_{\mu,j}\in\N_0$ satisfy the inequality
		\setcounter{equation}{4}
		\begin{align}\label{eq:GrahlCondition}
			(k-1)\cdot\sum\limits_{j=1}^{s_\mu}\omega_{\mu,j}+\ell\cdot s_\mu\leq \ell\cdot k
		\end{align}
		for all $\mu=1,\hdots,N$, where equality can only hold if $2\leq s_\mu\leq k-1$.\\
		Suppose that for every $f \in \F$ and every $z \in \D$, the following inequality holds:
		\begin{align}
			a(z)\cdot f^k(z)+f^{(\ell)}(z)+P[f](z)\neq b(z).
		\end{align}
		Then the family $\F$ is normal.
	\end{thmx}\medskip
	
	\begin{thmx}[\mbox{\cite[Theorem~4]{Grahl}}]\label{th:GrahlEntire}
		Let $f\in\mathcal{H}(\C)$, $\ell\in\N$ and $k\in\N\backslash\{1\}$. Let $a,a_1,\hdots,a_N\in\C$ with $a\neq 0$ and consider a differential polynomial 
		\begin{align}
			P[u]=\sum\limits_{\mu=1}^N a_\mu\cdot\prod\limits_{j=1}^{s_\mu} u^{(\omega_{\mu,j})},
		\end{align}
		where $2\leq s_\mu\leq k-1$ and $\,\sum_{j=1}^{s_\mu}\omega_{\mu,j}\neq 0$ holds for all $\mu=1,\hdots,N$.\\
		If the function $a\cdot f^k+f^{(\ell)}+P[f]$ does not vanish, then $f$ must be constant.
	\end{thmx}\medskip
	\noindent
	We conclude this section by proving the claim stated at the end of Section~\ref{se:Main}.\medskip
	
	\begin{prop}\label{pr:Counterexample}
		There exists $f\in\mathcal{M}(\C)$ such that $f$ is not of the form $f(z)= \mfrac{a e^{\alpha z} +b}{ce^{\alpha z} +d}$ for any choice of $a,b,c,d,\alpha\in\C$, while $S_2(f)$ omits the values $0$ and $\infty$.
	\end{prop}\medskip
	\noindent
	\begin{proof}
		We consider the Bessel functions defined by (cf. \cite[9.1.12 and 9.1.13]{AbramowitzStegun})
		\begin{align}
			J_0(z)&\coloneq\sum_{k=0}^\infty \frac{(-1)^k}{(k!)^2} \Big(\frac{z}{2}\Big)^{2k} \qquad\qquad\text{and}\\
			Y_0(z)&\coloneq\frac{2}{\pi}\big(\log\frac{z}{2}+\gamma\big)J_0(z)+\frac{2}{\pi}\sum_{k=1}^\infty\frac{(-1)^{k+1}H_k}{(k!)^2}\Big(\frac{z}{2}\Big)^{2k},
		\end{align}
		where $H_k$ denotes the $k$-th harmonic number and $\gamma$ is Euler's constant. These are linearly independent solutions of the Bessel differential equation (cf. \cite[9.1.1]{AbramowitzStegun})
		\begin{align}
			z^2y''(z)+zy'(z)+z^2y(z)=0.
		\end{align}	
		Next, we define
		\begin{align}
			f_1(z)\coloneq J_0(e^{z/2})\qquad\text{and}\qquad f_2(z)\coloneq Y_0(e^{z/2}).
		\end{align}
		Due to the logarithmic singularity of $Y_0$, we initially restrict our consideration of $f_2$ to the horizontal strips
		\begin{align}
			S_n\coloneq\{z\in\C\,:\,(4n-2)\pi<\text{Im}(z)<(4n+2)\pi\}.
		\end{align}
		Within each $S_n$, a direct computation yields
		\begin{align}
			f_2''(z)+\frac{e^z}{4}\,f_2(z)
			=\frac{e^z}{4}\,Y_0''(e^{z/2})+\frac{e^{z/2}}{4}\,Y_0'(e^{z/2})+\frac{e^z}{4}\,Y_0(e^{z/2})=0,
		\end{align}
		which implies that $f_2$ satisfies a second-order linear differential equation with entire coefficients. Hence, $f_2$ is analytically extensible to an entire function (\mbox{cf. \cite[Satz 3.2]{Herold}}).
		Similarly, $f_1$ satisfies 
		\begin{align}
			f_1''(z) + \frac{e^z}{4}\,f_1(z) = 0.
		\end{align}
		Now, a classical result about the Schwarzian derivative (see \cite[Theorem~6.1]{Laine}) implies that the function $f\coloneq f_1/f_2$ fulfills $S_2(f)(z)=e^z/2$. In particular, $S_2(f)$ omits the values $0$ and $\infty$.\\
		Next, we suppose that $f(z)=\frac{a e^{\alpha z} +b}{ce^{\alpha z} +d}$ for some choice of $a,b,c,d,\alpha\in\C$. Then the poles and zeros of $f$ must be periodic. However, let $j_{0,n}$ and $y_{0,n}$ denote the $n$-th positive real zero of $J_0$ and $Y_0$ respectively. Then asymptotically (cf. \cite[9.5.12]{AbramowitzStegun})
		\begin{align}
			j_{0,n}\sim\Big(n-\frac{1}{4}\Big)\pi \qquad\text{and}\qquad y_{0,n}\sim\Big(n-\frac{3}{4}\Big)\pi.
		\end{align}
		Thus, the positive real zeros of $J_0$ and $Y_0$ are almost equidistant, which shows that the zeros and poles of $f$ are not periodic. Hence, $f$ can not be expressed as a composition of a Möbius transformation and an exponential function.
	\end{proof}

	\section{Proofs of the Main Results}	
	\noindent
	Throughout this section, we consider $k\in\N$ and $M\in\R$ fixed. \medskip
	
	\noindent
	\textit{Proof of Theorem~\ref{th:HkM}.}
	\begin{ProofOf}\noindent
		From Theorem~\ref{th:ChuaquiGrönRättyä} we know that the derivative of every $f\in\HkM$ can be written as $f'=1/h^k$ for some zero-free function ${h\in\HD}$. Additionally we have ${S_k(f)=-k\,h^{(k)}/h}$. Therefore, the family ${\{h^{(k)}/h\in\HD\,:\,f'=1/h^k \text{ for some }f\in\HkM\}}$ is locally uniformly bounded. Applying Theorem~\ref{th:Schwick} shows that ${\{h'/h\in\HD\,:\,f'=1/h^k \text{ for some }f\in\HkM\}}$ is locally uniformly bounded as well. Since $f''/f'=-k\, h'/h$, it follows that the family $\HkMlog$ is locally uniformly bounded.\\
		Next, we consider a sequence $\fn\subseteq\HkM$. Because $(f_n''/f_n')_n$ is locally uniformly bounded and due to Proposition~\ref{pr:Poles}, we are able to apply Lemma~\ref{le:Lift} to a subsequence of $\fn$ with $E=\emptyset$. This implies that $\HkM$ is quasi-normal, that $\HkML$ is normal, and that no sequence in $\HkML$ converges to $\infty$. \pagebreak
		
		\noindent For the final claim, we use the inequality $1+x< 2\,(1+x^2)$ for $x\in\R$ to obtain
		\begin{align}\label{eq:Marty}
			(f')^\#(z)=\frac{|f''(z)|}{1+|f'(z)|^2}<\frac{2\,|f''(z)|}{1+|f'(z)|}<2\, \bigg\vert\frac{f''(z)}{f'(z)}\bigg\vert
		\end{align}
		for all $f\in\MD$ and $z\in\D$. Since $\HkMlog$ is locally uniformly bounded, it follows that the spherical derivatives of the functions in $\HkMD$ are locally uniformly bounded. Now, Marty's theorem implies that $\HkMD$ is normal.
	\end{ProofOf}\medskip
	
	\noindent
	\textit{Proof of Theorem~\ref{th:MkM}.}
	\begin{ProofOf}\noindent
		We consider a sequence $\fn\subseteq\MkM$. By Lemma~\ref{le:FinitePoles}, there exists $N\in\N$, such that each $f_n$ has at most $N$ poles. Therefore, we can find a subsequence $(f_{n_k})_{n_k}$ and a set $E\subseteq\D$ consisting of at most $N$ points, such that for every $z\in\D\backslash E$ there exists a neighborhood $U$ of $z$ in which almost all $f_n$ are analytic. Applying Theorem~\ref{th:HkM} locally on $\D\backslash E$ implies that there is a subsequence $(f_{n_\ell}''/f_{n_\ell}')_{n_\ell}\subseteq(f_{n_k}''/f_{n_k}')_{n_k}$ that converges locally uniformly on $\D\backslash E$ to some $F\in\mathcal{H}(\D\backslash E)$. Since $E$ has no accumulation point, we conclude that $\MkMlog$ is quasi-normal and that no sequence in $\MkMlog$ converges to $\infty$.\\
		Analogously to the proof of Theorem~\ref{th:HkM}, applying Lemma~\ref{le:Lift} to $(f_{n_\ell})_{n_\ell}$ shows that both $\MkM$ and $\MkML$ are quasi-normal, and that no subsequence of $(f_{n_k}'/f_{n_k})_{n_k}$ converges to $\infty$.\\
		Finally, the convergence of $(f_{n_\ell}''/f_{n_\ell}')_{n_\ell}$ on $\D\backslash E$ to $F\in\mathcal{H}(\D\backslash E)$ together with inequality~\eqref{eq:Marty} show that $((f_{n_\ell}')^\#)_{n_\ell}$ is locally uniformly bounded on $\D\backslash E$. Therefore, $\MkMD$ is quasi-normal by Marty's theorem.
 	\end{ProofOf}\medskip
 	
	\noindent
	\textit{Proof of Lemma~\ref{le:SchwarzianDifferentialPolynomial}.}
	\begin{ProofOf}\noindent
		Since $f$ is non-constant, there exists a neighborhood $U\subseteq\D$ such that ${S_k(f)\vert_U\in\mathcal{H}(U)}$. By Theorem~\ref{th:ChuaquiGrönRättyä}, there exists a zero-free ${h\in\mathcal{H}(U)}$ with $f'=1/h^k$ on $U$. Furthermore, we can find a logarithm $L\in\mathcal{H}(U)$ of $h$, i.e. $h=\e^L$. Using Faà di Bruno's formula (cf. \cite[Chapter 2.8]{Riordan}), we have that
		\begin{align}
			h^{(k)}=\big(\e^L\big)^{(k)}=\mathlarger{\sum}_{\mathclap{(n_1,\hdots,n_k)\,\in\,\Lambda}}\,\,\,\,\,\,\Bigg(\frac{k!}{n_1!\cdot\hdots\cdot n_k!} \,e^{L(z)}\,\prod\limits_{j=1}^k\bigg(\frac{L^{(j)}(z)}{j!}\bigg)^{n_j}\Bigg),
		\end{align}
		where $\Lambda$ is the set of all natural tuples $(n_1,\hdots,n_k)$ with $\sum_{j=1}^{k}j\cdot n_j=k$. Now, using Theorem~\ref{th:ChuaquiGrönRättyä} and $L'=h'/h=-f''/(k\cdot f'')$, we have
		\begin{align}\label{eq:SchwarzianDifferentialPolynomial}
			\begin{split}
				S_k(f)
				=-k\,\frac{h^{(k)}}{h}&=\mathlarger{\sum}_{\mathclap{(n_1,\hdots,n_k)\,\in\,\Lambda}}\,\,\,\,\,\,\Bigg(\frac{-k\cdot k!}{n_1!\cdot\hdots\cdot n_k!}\,\prod\limits_{j=1}^k\bigg(\frac{L^{(j)}(z)}{j!}\bigg)^{n_j}\Bigg)\\
				&=\mathlarger{\sum}_{\mathclap{(n_1,\hdots,n_k)\,\in\,\Lambda}}\,\,\,\,\,\,\Bigg(\frac{-k\cdot k!}{n_1!\cdot\hdots\cdot n_k!}\,\prod\limits_{j=1}^k\Bigg(\frac{1}{-k\cdot j!}\cdot\bigg(\frac{f''}{f'}\bigg)^{(j-1)}\Bigg)^{n_j}\Bigg)\\
				&=\mathlarger{\sum}_{\mathclap{(n_1,\hdots,n_k)\,\in\,\Lambda}}\,\,\,\,\,\,\Bigg(\frac{-k\cdot k!}{\prod_{j=1}^k \big((-k)\cdot (j!)\big)^{n_j}\cdot n_j!}\,\prod\limits_{j=1}^k\Bigg(\bigg(\frac{f''}{f'}\bigg)^{(j-1)}\Bigg)^{n_j}\Bigg).
			\end{split}
		\end{align}
		Finally, by the identity theorem, this equality extends from $U$ to $\D$.
	\end{ProofOf}\pagebreak
	
	\noindent
	\textit{Proof of Theorem~\ref{th:Exeption}.}
	\begin{ProofOf}\noindent
		Let $\Lambda$ be defined as in Lemma \ref{le:SchwarzianDifferentialPolynomial}, and set ${\tilde{\Lambda}\coloneq\Lambda\backslash\{(k,0,\hdots,0),(0,\hdots,0,1)\}}$. Next, we fix an enumeration $\mu$ of $\tilde{\Lambda}$, denote $N\coloneq \vert\tilde{\Lambda}\vert$ and define for each multi-index $(n_1,\hdots,n_k)\in\tilde{\Lambda}$ that
		\begin{align}
			a_{\mu^{-1}(n_1,\hdots,n_k)}&\coloneq \frac{-k\cdot k!}{\prod_{j=1}^k \big((-k)\cdot (j!)\big)^{n_j}\cdot n_j!}.
		\end{align}
		In addition, for $\nu=\mu^{-1}(n_1,\hdots, n_k)$ we set:
		\begin{align}
			s_\nu&\coloneq \sum_{r=1}^k n_r &&\text{and for} &&j=1,\hdots, s_\nu,\qquad\text{we define}\\
			\omega_{\nu,j}&\coloneq r-1  &&\text{whenever} &&n_1+\hdots, n_{r-1}<j\leq n_1+\hdots+n_r.
		\end{align} 
		Using Lemma~\ref{le:SchwarzianDifferentialPolynomial}, we get for all $f\in\F$ that
		\begin{align}
			S_k(f)
			&=\frac{(-1)^{k+1}}{k^{k-1}}\,\bigg(\frac{f''}{f'}\bigg)^k+\bigg(\frac{f''}{f'}\bigg)^{(k-1)}\!+\,\,\,\mathlarger{\sum}_{\mathclap{(n_1,\hdots,n_k)\,\in\,\tilde{\Lambda}}}\,\, a_{(n_1,\hdots,n_k)}\,\prod\limits_{j=1}^k\Bigg(\bigg(\frac{f''}{f'}\bigg)^{(j-1)}\Bigg)^{n_j}\\	\label{eq:FormGrahl}	
			&=\frac{(-1)^{k+1}}{k^{k-1}}\,\bigg(\frac{f''}{f'}\bigg)^k+\bigg(\frac{f''}{f'}\bigg)^{(k-1)}+\underbrace{\sum_{\mu=1}^N a_\mu\,\prod\limits_{j=1}^{s_\mu}\bigg(\frac{f''}{f'}\bigg)^{(\omega_{\mu,j})}}_{\eqcolon P[f''/f']}.
		\end{align}
		Now, we want to apply Theorem~\ref{th:GrahlNormal} with $\ell\coloneq k-1$. Note that $f''/f'$ is analytic, since $f$ is locally univalent. Furthermore, condition \eqref{eq:GrahlCondition} holds for $\mu=1,\hdots, N$, because $\sum_{r=1}^k n_r\cdot r=k$ for $(n_1,\hdots,n_k)\in\tilde{\Lambda}$ and therefore we get
		\begin{align}
			(k-1)\cdot\sum\limits_{j=1}^{s_\mu}\omega_{\mu,j}+\ell\cdot s_\mu 
			=(k-1)\sum\limits_{r=1}^{k}n_r\,(r-1)+(k-1)\,\sum_{r=1}^k n_r
			=\ell\cdot k
		\end{align}
		independently of $\mu$.
		Additionally, we have $2\leq s_\mu\leq k-1$ for all $\mu=1,\hdots, N$, because we removed the multi-indices $(k,0,\hdots,0)$ and $(0,\hdots,0,1)$ from $\tilde{\Lambda}$. Hence, Theorem \ref{th:GrahlNormal} applies and we conclude that $\F''/\F'$ is a normal family. 
	\end{ProofOf}\medskip
	
	\noindent
	\textit{Proof of Corollary~\ref{co:ExceptionalEnitre}.}
	\begin{ProofOf}\noindent
		We reuse the notation used in the proof of Theorem~\ref{th:Exeption} and write $S_k(f)$ as in equation~\eqref{eq:FormGrahl}. By Proposition~\ref{pr:Poles}, $f'$ does not vanish and therefore $f''/f'$ is an entire function. Since $2\leq s_\mu\leq k-1$ and $\,\sum_{j=1}^{s_\mu}\omega_{\mu,j}\neq 0$ holds for all $\mu=1,\hdots,N$, we are able to apply Theorem~\ref{th:GrahlEntire} to $f''/f'$. Thus, $f''/f'\equiv b$ for some constant $b\in\C$. Notice that $b\neq 0$, since $S_k(f)$ omits the value $0$. Therefore, $f$ must be of the form $f(z)=ae^{bz}+c$ for some $a,b,c\in\C$ with $a,b\neq 0$.
	\end{ProofOf}
		
	\section*{Acknowledgment}
	\noindent
	I am very grateful to Jürgen Grahl for several valuable remarks.

	\bibliographystyle{amsplain}

\end{document}